On Generating A Diminimal Set of Polyhedral Maps on the Torus

By

JENNIFER JOYCE HENRY
B.S. (University of San Francisco) 1995
M.A. (University of California, Davis) 1997

DISSERTATION

Submitted in partial satisfaction of the requirements for the degree of

DOCTOR OF PHILOSOPHY

in

Mathematics

in the

OFFICE OF GRADUATE STUDIES

of the

UNIVERSITY OF CALIFORNIA

DAVIS

Approved:

_________________________________

_________________________________

_________________________________

Committee in Charge

2001

<u>Abstract</u>

We develop a method to find a set of diminimal polyhedral maps on the torus from which all other polyhedral maps on the torus may be generated by face splitting and vertex splitting. We employ this method, though not to its completion, to find  53  diminimal polyhedral maps on the Torus.

## Introduction

Minimal polyhedral maps in manifolds have been given some attention. Steinitz found there to be only one on the sphere [13]. Barnette found there to be only seven on the projective plane [3]. For no other manifold has the entire set been found. Riskin has treated the pinched torus and began work on the torus [11 and 12]. In this paper, we take up work on the torus. We began by developing a set of maps from which all of the diminimal polyhedral maps may be generated. A computer is used to add edges in all necessary ways in order to arrive at a set of diminimal polyhedral maps.

## Preliminaries

We consider graphs embedded in the torus with no multiple edges, no loops, and no vertices of valence less than three. If $G$ is a graph embedded in a surface, then the closure of each connected component of the complement of the graph is called a *face*. An embedded graph together with its embedding is called a *map* if each vertex of the graph is at least 3-valent and each face is a closed 2-cell.

A *region* is a connected union of one or more faces. A *planar region* is a region that is simply connected. If the intersection of two faces is empty or connected then the faces are said to *meet properly*. If each face is simply connected and each pair of faces meet properly, then the map is a *polyhedral map*. We will refer to a polyhedral map on the torus as a *toroidal polyhedral map*, or TPM.

A graph embedded in the torus can be represented by drawing the graph on a rectangle whose top and bottom sides are identified, and whose left and right sides are identified. If an edge of the graph coincides with a section of the boundary of the rectangle, then in the figure, that section of the boundary will be solid, all other sections of the boundary of the rectangle in the figure will be dotted.  Such a representation of an embedded graph is called a *rectangular representation*.

The phrase *edge removing* refers to the process of obtaining one map from another by removing a single edge.  If the removal creates a 2-valent vertex then the two edges adjacent to that vertex are coalesced into one edge in the new map.  Let $G$ be a TPM and let $G'$ be the map obtained from $G$ by removing edge $e$ from $G$. If $G'$ is also a TPM, then edge $e$ is called *removable*.

The phrase *edge shrinking* refers to the process of obtaining one map from another by contracting an edge so that its two vertices become one vertex.  If the contraction creates a two-sided face in the new map then the two edges of the two-sided face are replaced with a single edge.  Let $G$ be a TPM and let $G'$ be the map obtained from $G$ by shrinking edge $e$ in $G$. If $G'$ is also a TPM, then edge $e$ is called *shrinkable*.

If $G$ is a TPM with no shrinkable edges and no removable edges then $G$ is called *diminimal*.  To generate diminimal TPM's we will use a process called edge adding (this process is also called face splitting in some papers).  Edge adding is the inverse of edge removing.  It is a process by which a map $G'$ is obtained from a map $G$ by adding an edge $e$ across a face of $G$ such that no edge of $G$ contains both vertices of $e$.

*Vertex splitting* is the inverse of edge shrinking.  It is a process by which a map $G'$ is constructed from a map $G$ by splitting one vertex, $x$, into two vertices, $x_1$ and $x_2$, and assigning each edge incident to $x$ to be incident to $x_1$, $x_2$, or both, and adding a new edge from $x_1$ to $x_2$. The process of vertex splitting is the dual to the process of edge adding. That is, if $G'$ is constructed from $G$ by adding an edge, the dual of $G'$ can be constructed from the dual of $G$ by splitting one vertex.

A *closed path* is an alternating sequence of vertices and edges,

$$v_1, e_2, v_3, \ . \ . \ . \ e_n, v_{n+1}$$

such that each $e_i$ is incident with $v_{i-1}$ and $v_{i+1}$, and such that no vertex appears more than once in the sequence with the exception that the first vertex may also be the last vertex in the sequence. An *open path* is an alternating sequence of vertices and edges, $e_2, v_3, \ . \ . \ . \ e_n$, such that each $e_i$ is incident with $v_{i-1}$ and $v_{i+1}$, and such that no vertex appears more than once in the sequence. When we denote a path we may use just the vertices, $v_1, v_3, \ . \ . \ . \ v_{n+1}$. Unless otherwise states, our paths will be open. If $U$ is $u_1$, $a_2, u_3, \ . \ . \ . \ a_n, u_{n+1}$ and $V$ is $v_1, b_2, v_3, \ . \ . \ . \ b_m, v_{m+1}$ and if vertex $u_{n+1}$ is vertex $v_1$, and no other $u_i$ is incident with any $v_j$, then a path $P$ may be written as $U \cup V$ and will be $u_1, a_2, u_3, \ . \ . \ . \ a_n, v_1, b_2, v_3, \ . \ . \ . \ b_m, v_{m+1}$. If two paths do not share more than their end vertices, then the paths are called *independent*.

A closed path with identical first and last vertex is called a *cycle*. The *length* of a path or cycle is the number of edges in its alternating sequence. If a cycle divides the torus into two regions, and one of those regions is a planar region, then the cycle is a *planar cycle*.

If a map $M$ can be obtained from a map $M'$ by replacing edges of $M'$ with paths then $M$ is said to *be a refinement* of $M'$. A map $M$ is said to *contain a refinement* of a map $M'$ if there is a refinement $J$ of the graph of $M'$ which is isomorphic to a subgraph $K$ of the graph of $M$ such that $J$ and $K$ have isomorphic embeddings.

**The Plan**

Let  T  be the set of all diminimal TPM's and let  S  be a set of maps such that every map in  T  contains a refinement of a map in  S.  If a map  M  contains a refinement of a map  M',  then  M  can be obtained from  M'  by edge adding  [see 2, Lemma 1]. Observe that since every map in  T  contains a refinement of a map in  S,  every map in  T may be generated by adding edges to the maps in  S.  Adding edges to the maps of  S  in all possible ways would generate not only the maps of set  T  but many others as well. We will develop a set of rules to determine each of the following, in order to shorten the process of generating  T  from  S.

- What edges must be added to the maps of   S  in order to generate the maps of   T?

- At any given step, which generated maps belong to  T?

- At any given step, which generated maps will not generate any map of  T?

Once the rules are established, a computer is employed to do the work of applying those rules in order to generate the diminimal TPM's.

**A Starting Set S**

We want to develop a set  S  such that every map in  T  contains a refinement of one of the maps in  S.  It has been proven by Barnette that every  TPM  is the union of two face disjoint annular regions, which we will call *bands*  [see 1, Theorem 4].  Riskin found all diminimal  TPM's  containing three bands and showed that there are no diminimal TPM's containing more than three bands  [12].  Therefore, let us redefine  T  to be the set of all diminimal  TPM's  containing precisely two bands.  When we have found the members of  T  and combined them with the members of Riskin's set we will have a complete set of diminimal  TPM's.

Let  G  be an arbitrary map in  T.  Map  G  will contain two bands.  Consider a rectangular representation of  G  like that shown in Figure  1.  Let us construct a graph  H on the sphere as follows.  Consider the rectangular representation of  G  as a planar graph.  Place two vertices  x  and  y  to the left and right respectively of the representation of  G.  Connect  x  to the vertices of the left hand path of  G  and  y  to the vertices of the right hand path of  G  as shown in Figure  2.  Next, identify the top of the picture (from x,  along the top of the representation of  G,  to  y)  with the bottom of the picture (from x,  along the bottom of the representation of  G,  to  y).  The result is a graph on the sphere.  Since  G  was a polyhedral map, the faces in the rectangular representation meet properly and clearly the faces that were added meeting  x  and  y  meet each other and the faces in the bands properly, therefore  H  is the graph of a polyhedral map on the sphere. Since  H  is polyhedral,  H  must be  3-connected  [see 12,  Theorem 4]  and so there must be three independent paths from x  to  y.  That means across the original rectangular representation of  G  there must be three disjoint paths from left to right.  In particular this provides three disjoint paths across each of the two bands of  G.  We define  S  to be the set of all toroidal maps containing precisely two bands and three disjoint edges across each band and nothing more.  Every map of   T  will contain a refinement of one of the maps of  S.

To build  S  we begin with two bands.  Across one band, we put three paths (Figure  3).  Next, we add three paths across the remaining band in all possible ways. The first path across the remaining band may lead from an existing vertex to another existing vertex, from an existing vertex to an existing edge, or from an existing edge to another existing edge. These three types of paths will be referred to as vv, ve, and ee respectively.

Upon adding the first edge to the second band, there are two possible vv paths up to isomorphism. Similarly, there are two possible ev paths and two possible ee paths. These maps are shown in Figure 4 labeled A through F and will be referred to as such.

We now consider the choices for adding a second edge to the second band of A. We add edges in all possible ways from the path (1,2,3,1) to the path (4,5,6,4). The computer does this job and provides us with 25 new maps.

Let G and H be two isomorphic maps generated by this process. Let $e$ be the first edge added to the second band of G. Let $v_1$ be the vertex of $e$ on cycle (1,2,3) and $v_2$ be the vertex of $e$ on cycle (4,5,6).

**Definition**: If R is the face in G that borders $e$ in such a way that when traversed counterclockwise, vertex $v_1$ follows vertex $v_2$, then R is called the *rear-face*. If F is the face in G that borders $e$ in such a way that when traversed counterclockwise, vertex $v_2$ follows $v_1$ then F is called the *fore-face*.

**Definition**: If the isomorphism from G to H carries the fore-face of G to the fore-face of H, then we call the maps *simply-isomorphic*. If the isomorphism from G to H carries the fore-face of G to the rear-face of H, then we call the maps *complexly-isomorphic*. See Figure 5 for an example of each.

If two maps are found to be simply-isomorphic, one is removed from the set. If two maps are found to be complexly-isomorphic, both are left in the set. The reason is that in the next phase a third edge will be added to this band in the fore-face.

In addition to removing extraneous simply-isomorphic maps we also remove those containing the configurations described in Figure 6. It is easily seen that when a third edge is added to the fore-face of those configurations the resulting maps will contain three bands, see Figure 7 for some examples. Any diminimal maps constructed from these would also contain three bands and thus would have been found by Riskin.

To each of the maps of this shortened set the computer adds edges in all possible ways from the path $(x \ldots 1)$ to the path $(y \ldots 4)$ in the fore-face where $x$ and $y$ are the end vertices of the second edge. From the resulting set, extraneous isomorphic maps and all maps with three bands are removed. The resulting set of 18 maps make up the first portion of the $S$.

**Definition**: A *straight-across edge* is a $vv$ edge across the second band from vertex $j$ to vertex $j+3$ where $j \in (1,2,3)$.

Next, consider map B. Any map of $S$ that has a straight-across edge was generated (up to isomorphism) by map A. Therefore, as edges are added to maps B through F, no straight-across edges are allowed. Aside from that, the process for map B is the same as that of map A. The resulting 37 maps are added to the growing set $S$.

Any map containing a $vv$ edge across the second band was generated by map A or B, therefore we add no $vv$ edges to maps C through F. For each map C through F the computer is used to add a second edge from path $(1,2,3,1)$ to path $(7,5,6,4,7)$, and a third edge from path $(8,\ldots1)$ to path $(9,\ldots7)$, in all possible ways, excluding those ways that are covered by any earlier map. From this set extraneous isomorphic maps and maps with three bands are removed. The resulting maps are added to set $S$, which finally contains 359 maps.

**Lemma 1**. *Let  G  be an arbitrary member of  T.  Then graph  G  contains a band  B bordered by two cycles,  X  and  Y,  with three disjoint edges,  e₁,  e₂,  and  e₃,  such that each edge lies across  B  with one vertex on  X  and the other on  Y.*

**Proof**: Assume the lemma is false and let  G  be a member of  T  that does not satisfy the lemma.  We have already established that each member of  T  contains two bands, let us refer to them as  A  and  B.  We have also established that across each band there are three disjoint paths.  Select two bands in  G  so that the combined length of the paths across  B  is minimized.  Since  G  does not satisfy the lemma there must be a path  P  in G  across  B  of length at least two (see Figure  8).  If there exists a path  Q  from  some vertex  x  on the path  P  to a boundary vertex  y  of  B,  with  y  not on  P,  then  B  can be redrawn so that  x  is an endpoint of  P  (see Figure  9).  That would mean that we had not minimized the combined lengths of the paths across  B.  On the other hand, if there is no such path  Q  then  G  could be disconnected by removing the endpoints of  P  which is not possible since polyhedral maps are  3-connected  [see 11, Theorem 4].  ∎

A band which satisfies Lemma 1 and is minimal with respect to number of faces will be referred to as a *minimal band*.  We assume that paths  14,  25,  and  36 lie in a minimal band, and therefore will be edges in the final maps.

## Map Representation

It is important to arrive at a feasible representation of a map in order to write computer programs to manipulate the maps.  Although a graph is by definition a collection of vertices and a collection of edges (pairs of vertices), since we are interested in maps we also have to keep track of the embedding. Our solution is a variation on the well-known Rotational Embedding Scheme [7]. When using the Rotational Embedding Scheme, an embedding is given by specifying the order in which neighboring vertices surround each vertex. Whereas the Rotational Embedding Scheme relies on relationships among vertices, the method we develop relies on the relationship among vertices and faces. Each vertex in a map is assigned a number (beginning with one and incrementing.) Each face is recorded by listing in a counterclockwise direction the vertices that border it. The faces are put together into a one row matrix and each face is separated from the next by a -1. The first element of the matrix is the order of the underlying graph, the second is a zero (a placeholder) and the faces follow. For example, the map in Figure 10 might be represented as follows.

```
[14 0 -1 1 4 5 2 9 7 -1 2 5 8 6 3 11 -1 3 6 10 12
4 1 -1 14 13 8 5 4 12 -1 13 14 11 3 1 7 -1 13 7 9
10 6 8 -1 10 9 11 14 12 -1].
```

This method captures the same information captured in the Rotational Embedding Scheme and therefore is a feasible method for recording a map. Further, given the matrix for a particular map, one could easily obtain the Rotational Embedding Scheme for the dual of that map.

Matlab was used because it has many tools for manipulating matrices. Here we provide a general description of some of the operations performed to give the reader a sense of the style of our Matlab scripts. For this discussion, let $M$ be a matrix representation of a map. A cell array is used as a look up table in order that certain data be calculated once for each map. For example, the $i^{th}$ entry of the adjacency cell array is a vector containing the neighbors of vertex $i$. Often an operation must be performed for each face of a map. For that purpose an index vector is calculated using a find command on the number –1. With this index vector the i[th] face can quickly be found in M. A nested for loop is used when faces must be considered in pairs. For example, to determine if a map is polyhedral we consider how each pair of faces meet. The intersect command and size command are used on a pair of faces to get a general idea of how they meet. If they meet at exactly two vertices then the two faces are examined more closely to determine if the vertices form an edge in both faces. The map represented by M is polyhedral if and only if no two faces are found that meet at more than two vertices, or that meet at two vertices which do not form an edge in one of the two faces.

## Adding Edges

We have established a set S such that every map in T may be generated from at least one of the maps in S by edge adding. We examine the structure of the maps in T, specifically what being polyhedral implies, to determine the best way to add edges to the maps of S. We must add every edge that might produce a map of T, but we would also like to avoid doing more processing than is necessary.

**Lemma 2.** *Let G be a graph embedded in a 2-manifold M such that its faces are closed cells, and let E be a region of M whose boundary is $P_1 \cup P_2 \cup Q$, where $P_1$ and $P_2$ are edge disjoint paths and Q is either a vertex belonging to $P_1$ and $P_2$, or a path that is edge disjoint from $P_1$ and $P_2$. Let $P_1$ and $P_2$ meet at a vertex c, $P_1$ and Q meet at a vertex b and $P_2$ and Q meet at a vertex b' (where it might be that b = b'). If no face of G meets both Q and c then there is a path across E meeting the boundary of E only at its endpoints e and f with e on $P_1(b,c)$ and f on $P_2(b',c)$.*

This lemma was proven by Barnette [1] and will be useful for the next lemma.

**Lemma 3**. *Let $G$ be a map in $T$. Let $E$ and $F$ be face disjoint regions of $G$ such that the boundary of $E$ is $b \cup P_1 \cup c \cup P_2$, where $P_1$ and $P_2$ are edge disjoint paths meeting at vertices $b$ and $c$, and the boundary of $F$ is $b \cup Q_1 \cup c \cup Q_2$, where $Q_1$ and $Q_2$ are edge disjoint paths meeting at vertices $b$ and $c$. Then in $G$ either there is a path across $E$ meeting the boundary of $E$ only at its endpoints $e_1$ and $e_2$ with $e_1$ on $P_1(b,c)$ and $e_2$ on $P_2(c,b)$, or there is a path across $F$ meeting the boundary of $F$ only at its endpoints $f_1$ and $f_2$ with $f_1$ on $Q_1(b,c)$ and $f_2$ on $Q_2(c,b)$.*

**Proof**: Since $G$ is a polyhedral map, none of its faces meet improperly. Thus, if there is a face across $E$ meeting $b$ and $c$ then there is no face across $F$ meeting $b$ and $c$ and by Lemma 2 there is a path in $F$ as described. If there is no face across $E$ meeting $b$ and $c$ then there is a path in $E$ as described. ∎

**Lemma 4**. *Let $H$ be a map in the data set. Let $G$ be a map in $T$ containing a refinement of $H$. Let $E$ and $F$ be two faces of $H$. Let $E$ have boundary $Q \cup P_1 \cup c \cup P_2$, where $P_1$ and $P_2$ are edge disjoint paths and $Q$ is a path edge disjoint from $P_1$ and $P_2$. Let $P_1$ and $P_2$ meet at a vertex $c$, $P_1$ and $Q$ meet at a vertex $a$ and $P_2$ and $Q$ meet at a vertex $b$. Let $F$ have boundary $Q \cup R_1 \cup c \cup R_2$, where $R_1$ and $R_2$ are edge disjoint paths and $Q$ is a path edge disjoint from $R_1$ and $R_2$. Let $R_1$ and $R_2$ meet at $c$, $R_1$ and $Q$ meet at a vertex $b$ and $R_2$ and $Q$ meet at a vertex $a$. Then there is a refinement of $H$ in $G$ with faces $E^*$ and $F^*$ (corresponding to $E$ and $F$ under homeomorphism as shown in Figure 11) such that one of the ten configurations described in Figure 7 must be present.*

**Proof**: Let $EF$ denote the region in $H$ which is the union of faces $E$ and $F$. Construct a graph $K$ on the sphere as follows; consider a rectangular representation of $EF$ as a planar graph. Place two vertices $x$ and $y$ above and below the representation of $EF$ respectively. Connect $x$ to the vertices of path $c \cup R_2 \cup a \cup P_1 \cup c$ and $y$ to the vertices of the path $c \cup P_2 \cup b \cup R_1 \cup c$ as shown in Figure 12. Identify the left-hand side of the picture with the right hand side of the picture. The result is a 3-connected graph on the sphere (same argument used in the section *A Starting Set S*) so there must be three independent paths from $x$ to $y$. Only one of those paths can pass through vertex $c$, therefore there must be two disjoint paths across EF. Let us call these two paths U and V.

We proceed by considering how U and V meet Q. If at least one of U or V does not meet Q, then E* and F* are precisely E and F and one of the configurations A or B of Figure 7 is present. Otherwise, both U and V meet Q at least once. Construct a new path Q* as follows (for convenience we shall refer to one end of Q* as the head of Q* and the other end as the tail). Begin the head of Q* at vertex a and let it travel along Q until it reaches one of the paths U or V. Begin the tail of Q* at vertex b and let it travel along Q until it reaches one of the paths U or V. One of the following must occur.

*Case 1:* The head and tail of Q* have met the same path. Join the head and tail of Q* by traveling along that path. Faces E* and F* are formed by dividing EF along Q*. Since U and V do not meet each other, and Q* meets only one of them and one of the configurations A or B of Figure 7 is present.

*Case 2:* The head and tail of Q* have met different paths. Without loss of generality we will say that the head of Q* has met U at vertex u and the tail of Q* has met V at vertex v. Since U and V do not meet each other, there must be a segment of Q that passes from U to V. Let us denote this segment as $q_u q_v$ where $q_u$ is the vertex on U and $q_v$ is the vertex on V. Continue the head of Q* along U until reaching $q_u$. Continue the tail of Q* along V until reaching $q_v$. Join the head and tail of Q* with the segment $q_u q_v$. Faces E* and F* are formed by dividing EF along Q*. If vertex a is $q_u$ and vertex b is $q_v$ then one of the configurations C or D of Figure 7 is present. If vertex a is $q_u$ and vertex b is not $q_v$ then one of the configurations E or G of Figure 7 is present. If vertex a is not $q_u$ and vertex b is $q_v$ then one of the configurations F or H of Figure 7 is present. Finally if a is not $q_u$ and b is not $q_v$ then one of the configurations I or J of Figure 7 is present. ∎

Let $G_1$ be a member of S and $G_k$ be a member of T such that $G_k$ contains a refinement of $G_1$. Let $G_1$, $G_2$, . . . $G_k$ be a sequence of maps such that $G_{i+1}$ is constructed from $G_i$ by adding one or two edges across one or two faces of $G_i$. Since $G_k$ is a diminimal TPM, each $G_i$ in the sequence (i < k) is not a TPM. Therefore given $G_i$, two faces E and F can be found that meet improperly. Those two faces correspond to two regions E' and F' in $G_k$. Applying Lemmas 3 and 4 to $G_k$ will tell us that one of several path configurations must exist in $G_k$. For each possible path configuration in $G_k$ we generate a map $G_{i+1\_x}$ which contains an edge or pair of edges corresponding to that configuration. With such an edge adding process $G_i$ will generate one or more maps $G_{i+1\_a}$, $G_{i+1\_b}$ . . . $G_{i+1\_x}$. Each of these new maps is added to our data set for further consideration. The phrase *data set* will refer to the dynamic set of maps being considered. The set will begin as S, will increase in size as each element is replaced with all maps constructed from it ($G_{i+1\_a}$, $G_{i+1\_b}$ . . . $G_{i+1\_x}$), and will decrease in size as maps are removed from the set based on the rules established for recognizing members of T, maps that will not generate any member of T, and redundant isomorphisms.

## Maps To Throw Away

**Definition:** Let G be a member of our data set. Let C be a planar cycle in G. We say that G contains an *edge inside a cycle* (*e-i-c*) if there exists an edge *e* of G that is disjoint from C and lies in the planar region enclosed by C.

**Lemma 5**. *Let  $G_1$  and  $G_2$  be maps in our data set such that  $G_2$  may be generated from  $G_1$  by adding one edge across a face of  $G_1$. If  $G_1$  contains an  e-i-c  then so does  $G_2$.*

**Proof**:  Let  $G_1$  and  $G_2$  be as described in the hypothesis.  Let  C  be a planar cycle in  $G_1$  containing edge  *e*.  When an edge is added across a face of  $G_1$,  C  remains a planar cycle since the new edge will either not touch C, or will add one or two vertices to C thereby increasing its length and leaving it intact.  Similarly for  *e*,  the new edge will either not touch  *e*, or will add a vertex to  *e*  making it a path with both edges inside  C.  The result is that C remains a planar cycle with an edge inside.  ■

**Lemma 6.** *Let  $G_1$  and  $G_2$  be maps in our data set such that  $G_2$  may be generated from  $G_1$  by adding one edge across a face of  $G_1$. If the dual of  $G_1$  contains an e-i-c then so does the dual of  $G_2$.*

**Proof:**  Let  $G_1$  and  $G_2$  be as described in the hypothesis.  Let  $H_1$  and  $H_2$  be the dual maps respectively.  Since  $G_2$  can be constructed from  $G_1$  by adding one edge,  $H_2$  can be constructed from  $H_1$  by splitting one vertex.  Let  C  be a planar cycle in  $H_1$  containing edge  *e*.  If the vertex being split does not lie on  C,  then  C  remains a planar cycle.  If the vertex being split does lie on  C  then  C  still remains a planar cycle though its length may increase by one (see Figure  13).  Similarly for  *e*,  the new vertex will either leave  *e*  intact or split it into two edges, both of which are inside the cycle  C.  The result is that  C  remains a planar cycle with an edge inside.  ■

**Lemma 7.** *If  v  is a  3-valent vertex of a triangular face  F  of a  TPM  C  then the edge e  of  F  missing  v  is removable.*

This lemma was proven by Barnette for polyhedral maps on the projective plane [see 3, Lemma 3]. His proof does not rely on the topology of the projective plane and thus applies also to the torus. We will use this result in the following lemma.

**Lemma 8.** *If $G$ contains an e-i-c, then $G$ is not a diminimal TPM.*

**Proof**: Assume the statement is false and let $G$ be a diminimal TPM with an e-i-c. Among all cycles in $G$ in which an edge lies, let cycle $C$ with edge $e$ be minimal with respect to the number of faces in the planar region enclosed by $C$. Since $G$ is diminimal, $e$ is not shrinkable. Therefore, there must be two faces, A and B, that would meet improperly if $e$ were shrunk. Faces A and B will each meet distinct ends of $e$ and will meet each other at $f$, where $f$ is either an edge or a vertex. Let $v_a$ be the vertex where $e$ meets A, and let $v_b$ be the vertex where $e$ meets B. Let D be the planar region adjacent to A, B, and $e$ (see Figure 14). Let $f_d$ be used to denote $f$ if it is a vertex, and to denote the vertex of $f$ that meets region D if $f$ is a path. Let C' be the planar cycle that encloses A, B, D, and nothing else. Region D is not a triangular face since A and B were chosen to meet improperly when $e$ is shrunk.

Case 1: There is a vertex $g$ on path $v_a f_d$ or path $v_b f_d$, without loss of generality we use $v_a f_d$. There is not a second vertex, $h$, on $v_a f_d$ since such a vertex would provide an edge $gh$ inside the cycle C' contradicting our minimal cycle assumption. The vertex $g$ is not 2-valent in $G$ so there must be a third edge $gh$ in $G$. If $h$ lies inside D, or lies on $v_b f_d$ then $gh$ is an edge inside cycle C' contradicting our minimal cycle assumption. If $h$ is $v_a$ or $f_d$ we get a two sided face in $G$ which is not possible in a polyhedral graph (see Figure 15). The last choice for $h$ is that it be $v_b$, however the edge $gv_b$ would form the triangle $gv_bv_a$ with a three valent vertex $g$ (all other possible edges from $g$ have been ruled out) and by Lemma 7 that would make $e$ removable.

**Case 2:** There is no vertex on $v_af_d$ or $v_bf_d$. Since region D cannot be a triangular face there must be a vertex $g$ inside D. Vertex $g$ must be 3-valent and therefore would form the triangle $gv_bv_a$ and by Lemma 7 that would make $e$ removable.

Thus G cannot be a diminimal TPM. ∎

**Lemma 9.** *If the dual of G contains an e-i-c, then G is not a diminimal TPM.*

**Proof**: The dual of a diminimal TPM is a diminimal TPM. If the dual of G is not a diminimal TPM because it contains an e-i-c, then G is also not a diminimal TPM. ∎

**Corollary 10.** *Let G be a member of our data set. If G or its dual contain an e-i-c then no diminimal TPM can be generated from G by edge adding.*

**Proof**: If G contains an e-i-c then by Lemma 5 every map generated from G by edge adding will also contain an e-i-c, and by Lemma 8 no such map will ever be a diminimal TPM. If the dual of G contains an e-i-c then by Lemma 6 the dual of any map generated from G by edge adding will also contain an e-i-c and by Lemma 9 no such map will ever be a diminimal TPM. ∎

When a map G is found that has an e-i-c or such that the dual has an e-i-c the map is removed from the data set since it will not generate any map of T. When a map is found that has three bands it is removed from the data set since Riskin has generated all TPM's with three bands. We chose edges 14, 25, and 36 to be edges in our final maps so we also remove any map that does not contain those three edges. Lastly if $G_i$ and $G_j$ are isomorphic maps then one of them is removed from the data set.

## Minimal Polyhedral Maps

As we proceed to generate maps one important step is determining whether a map is polyhedral. To accomplish this we assume a map is polyhedral and look for a contradiction. We consider the faces in pairs. If a pair meets at fewer than two vertices then the pair meets properly. If a pair meets at two vertices and if these two are adjacent in both faces then the meeting is proper. Otherwise the meeting is obviously improper. If any improper meetings are found then the map is not polyhedral. If all meetings are proper then the map is polyhedral.

When a polyhedral map is found we next want to know if it is diminimal. We look for removable edges and shrinkable edges. We look for removable edges by considering the faces in pairs. If two faces meet along an edge then we consider the two faces as a single larger face and check to see if any of the other faces meet this larger face improperly. If no improper meetings are found, then this edge is removable and hence the map is not diminimal. If no removable edges are found we next look for shrinkable edges. To look for shrinkable edges we consider the edges one at a time. The edge $xy$ is shrunk by removing the vertex $y$ from the two faces containing the edge $xy$, and replacing $y$ with $x$ in any other face containing $y$. If the resulting map is polyhedral then the edge is shrinkable and the map is not diminimal. If no shrinkable edges are found then the map is diminimal.

When a map is found to be polyhedral and diminimal it is removed from the data set and added to our list of maps belonging to set  T.  When a map is polyhedral but not diminimal then it is removed from the data set since by definition no amount of edge adding will return it to diminimality.  When a map is found that is not polyhedral, then it must be assumed that adding more edges to it will generate a polyhedral map.  Such a map is removed from the data set and replaced by one or more maps generated from it by edge adding.  Each step removes an improper meeting without introducing any new improper meetings.  When all improper meetings have been removed, all maps are polyhedral, the data set is empty, and the process terminates.

## Isomorphism

It often happens that we arrive at two maps that are isomorphic to one another. We are not interested in keeping such maps so we develop a method for identifying a map that is isomorphic to another in order to remove it from the list.

Maps are considered in pairs  $G_i$  and  $G_j$.  The basic strategy is to look for properties that cause these maps to not be isomorphic.  If any of the following properties are not identical in such maps then they cannot be isomorphic;

- number of vertices

- number of faces

- number of  $k$  valent vertices (for any positive integer  $k$)

- number of  $k$  sided faces (for any positive integer  $k$)

In order for  $G_i$  to be isomorphic to  $G_j$  there must be a bijective mapping,  $f$,  from the set of vertices of  $G_i$  to the set of vertices of  $G_j$  such that vertices  $x$  and  $y$  are adjacent in $G_i$ if and only if vertices $f(x)$ and $f(y)$ are adjacent in $G_j$, and two faces in $G_i$ are adjacent if and only if their images are adjacent in $G_j$. If the isomorphism of two maps can not be ruled out by these facts then all feasible mappings from $G_i$ to $G_j$ are checked to see if they preserve all adjacencies. By feasible mappings we mean mappings that preserve vertex valence, face size, and a variety of adjacency constraints. If an isomorphism is found then we remove $G_j$ from the data set.

Here we would like to acknowledge that the method we have employed for finding isomorphisms might possibly not the most efficient for these particular maps. Much work has been done on the subject. Eugene Luks has done work regarding the isomorphism of graphs of bounded valence as well as parallel algorithems for graph isomorphisms [9 and 10].

## Results

The data set began with 359 maps. Over 12 months more than 800,000 maps have been processed and the data set has grown to over 1,000,000. Presented here are the diminimal maps we have found, statistical analysis regarding the complete set, and some conjectures.

Riskin found there to be only two diminimal TPMs containing three bands. Those two maps are shown in Figure 16. Riskin also found nine of the diminimal TPMs containing two bands. At this point, we have found 53 diminimal TPMs containing two bands. Those maps are described in Table 17.

When it became clear that we would not be able to complete the process in time for this paper, we decided to perform some statistical analysis on the data. Random samples of various sizes were selected from the data set and processed to completion. Table 18 contains a description of those samples. Although nine diminimal maps were generated, eight of those were isomorphic to maps we had already found. We believe that the diminimal set will have fewer than 60 maps. Barnette had conjectured in a private communication that there would be between 30 and 40 diminimal TPMs. Further, based on what we have seen with the maps processed to completion, we believe that no diminial map will have more than 17 vertices or 17 faces. Riskin did prove that there are a finite number of diminimal TPMs [11] but by no means does his proof indicate that the number is small.

Of the maps that we found there are seven that are self dual, see Figure 19. Each diminimal TPM has a dual TPM. We have indicated dual pairs as well as self dual maps in Table 20.

On the plane it is well established that the following two statements hold [14].

- If a graph has a polyhedral embedding in the surface, then every embedding in the surface is isomorphic to that embedding.

- If a graph has a polyhedral embedding in the surface, then a cycle is a face in one embedding if and only if it is a face in all other embeddings.

On the torus it is easily seen that the second statement does not hold by looking at $K_7$. We can now say that the first statement does not hold either.

**Theorem 11.** *There exists a graph that has two non-isomorphic polyhedral embeddings in the torus.*

**Proof**: Pictured in Figure 21 are two TPMs that are not isomorphic, but whose graphs are isomorphic. ■

When an edge is added to a map  G  of our data set, the number of vertices may go unchanged but the number of faces must increase.  The data set is sorted according to the number of faces in each map.  Those maps with the fewest faces were processed first.  At the time the processing was halted we had completed work on all maps with fewer than ten faces.  That provides us with the following.

**Theorem  12.**  *All diminimal TPMs with fewer than ten faces have been found.*

**Corollary  13.**  *All diminimal TPMs with fewer than ten vertices have been found.*

One of the criteria used to discard maps that we mentioned above was the criterion that there must be three paths across the rectangular representation of a map.  This criterion was developed after the process began.  If it had been applied at the beginning, the original set could have been decreased from  359  maps to  319  maps.  It is unclear how much processing this decrease would have allowed, but it should be considered when this effort is started again.

There is another method that we did not investigate soon enough to apply at the beginning, but should be considered when this effort is started again.  All the diminimal TPMs containing a pinched annulus can be identified with methods similar to those we have used.  With these identified, as the diminimal TPMs are generated the maps with a pinched annulus could be disregarded.

Here are some questions that can be taken up at a future date.

- How many of the diminimal TPMs are self dual?

- What is the upper bound on the number of vertices in a diminimal TPM?

- Do all diminimal TPMs contain a refinement of $K_{(3,3)}$ as shown in Figure 22? An answer of yes would imply that every TPM contains a refinement of this embedding. The author has made significant progress (in a separate work) in showing that every TPM contains such a refinement, and that there is no other closed 2-cell embedding on the torus such that every TPM contains a refinement of it. The argument is extremely long and tedious. Further, we have not found an efficient method for determining whether one map contains a refinement of another map.

- A TPM is called geometrically realizable if it can be embedded in Euclidean 3-space in such a way that each face is a flat polygon, no pair of neighboring faces lie in the same plane, and the resulting domain of $\mathbf{R}^3$ is bounded by finitely many hyperplanes supporting the polygonal facets. A good deal of attention has been given to geometrically realizable maps in general, and some have been found on the torus [6]. Which, if any, of the diminimal TPMs are geometrically realizable?

**Figures**

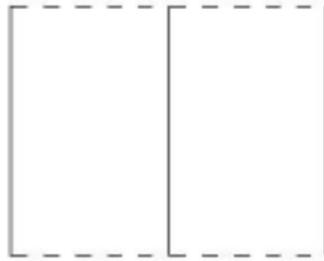

Figure 1  Rectangular representation of an arbitrary member of  T showing two bands

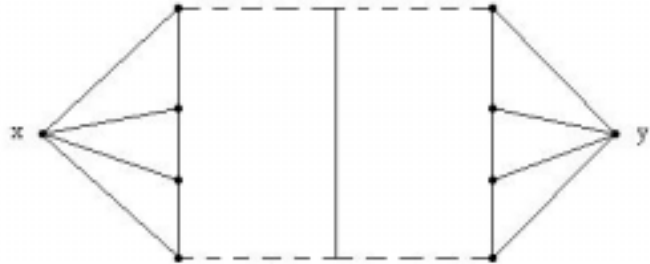

Figure  2  Planar map constructed from
rectangular representation of arbitrary member of T

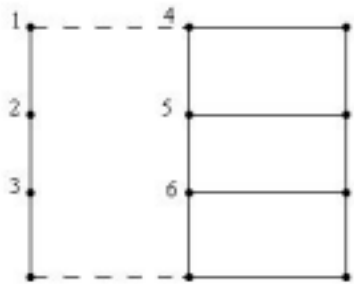

Figure  3  Three paths accross one band.

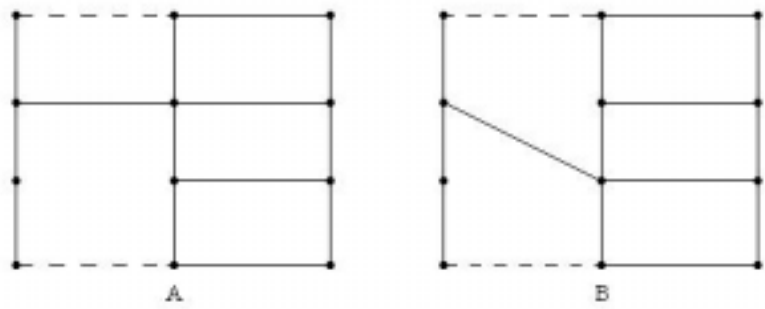

vv edges

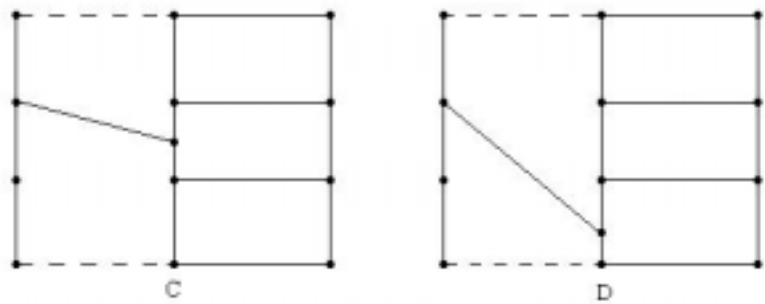

ve edges

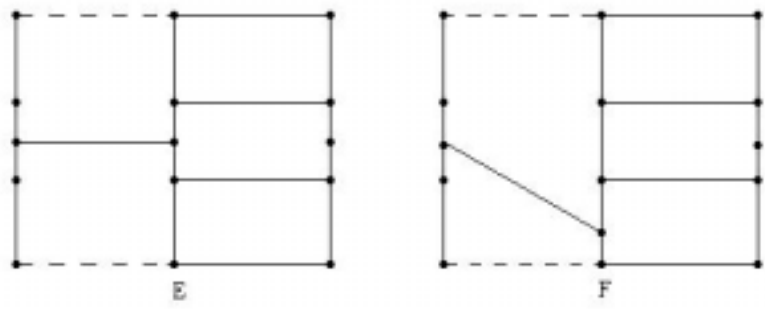

ee edges

Figure  4  Maps  A  through  F

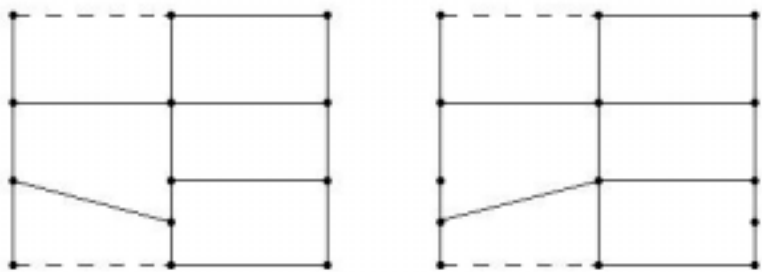

Simply-isomorphic

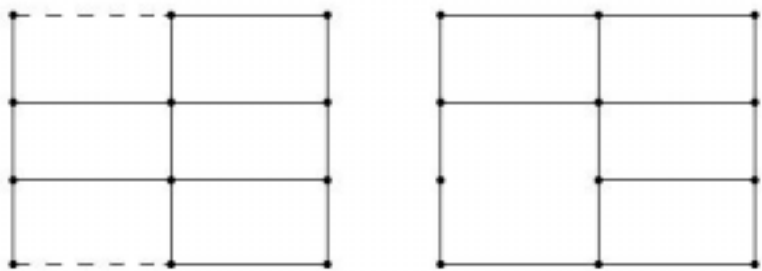

Complexly-isomorphic

Figure 5 Isomorphism examples

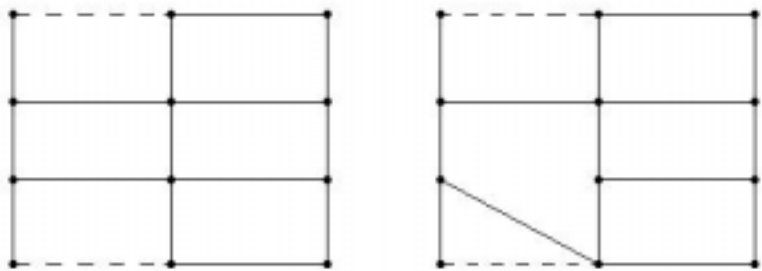

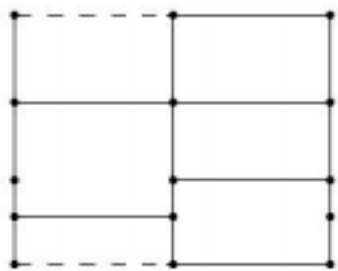

Figure 6 Three maps removed

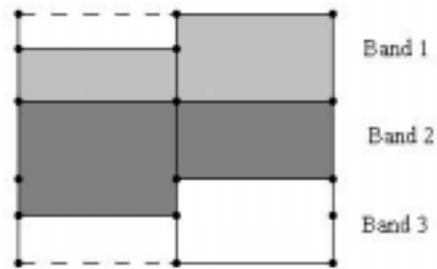

Figure  7  Example demonstrating three bands in generated map

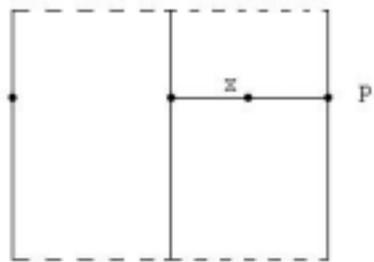

Figure  8  Path  P  accross a band with a vertex  x

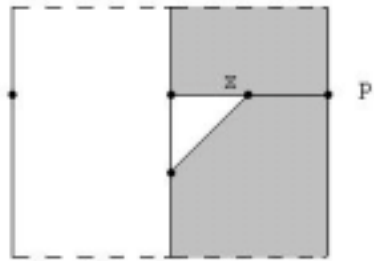

Figure  9  A path from x to the edge of the band
allows a new band to be found with fewer faces

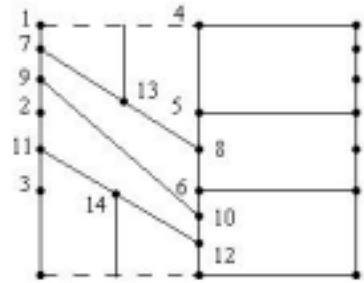

Figure 10  Example

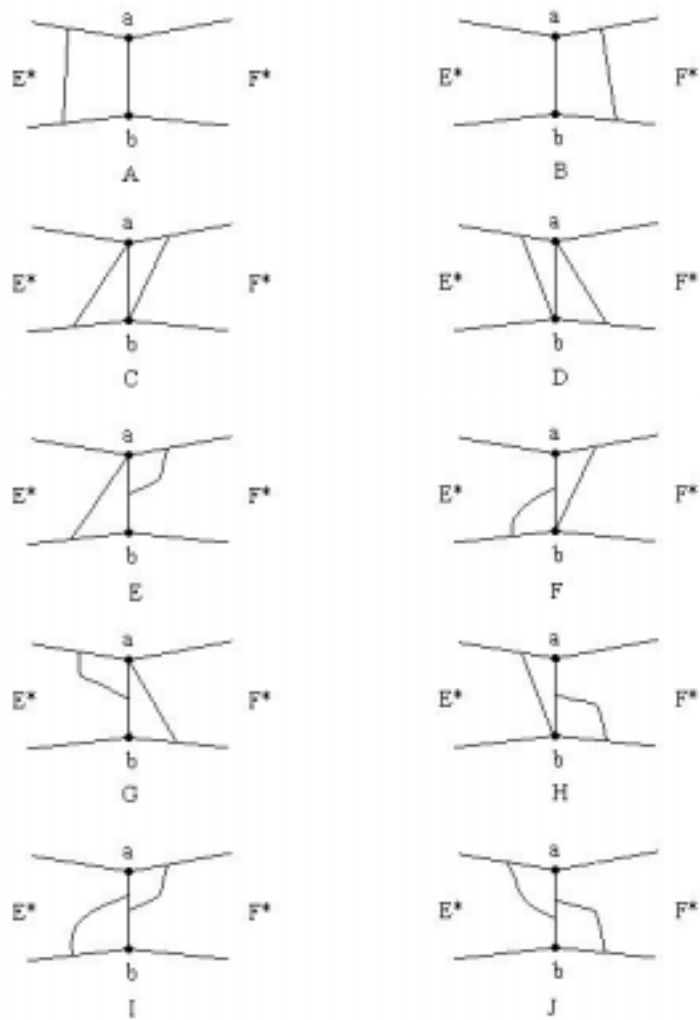

Figure 11  One of ten configurations must be present

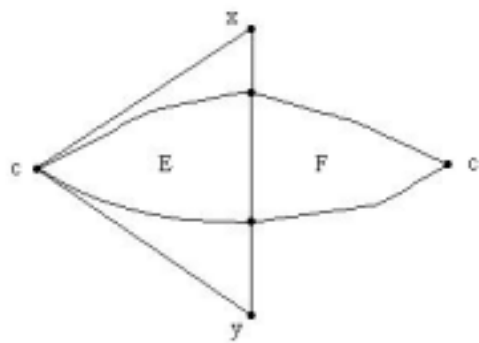

Figure 12

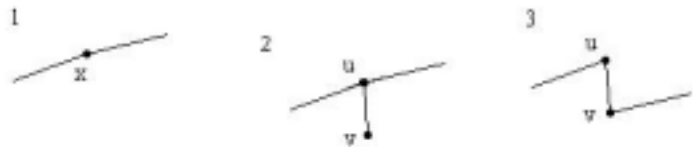

Figure 13  Vertex  x  lies on cycle C.
When  x  is split into  u  and  v,  C  is still intact

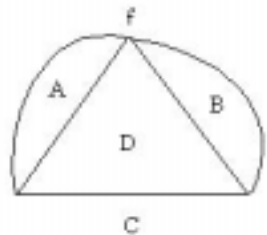

Figure 14  Depiction of region

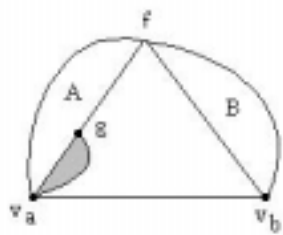

Figure 15  Vertex  g  causes a contradiction

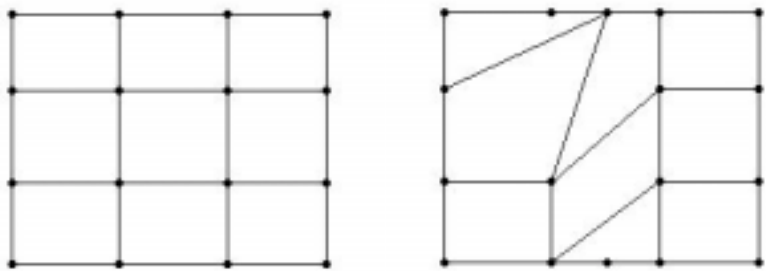

Figure  16  Riskin's diminimal TPMs with three bands

| # | Order | Size | P | V | Map |
|---|---|---|---|---|---|
| 1 | 7 | 14 | (14, 0) | (0, 0, 0, 7) | [ 7, 0, -1,  1,  4,  5, -1,  1,  5,  2, -1,  2,  5,  3, -1,  5,  6,  3, -1,  3,  6,  1, -1,  6,  4,  1, -1,  4, 7,  5, -1,  5,  7,  6, -1,  7,  1,  2, -1,  6,  7,  2, -1,  6,  2,  4, -1,  4,  2,  3, -1,  4,  3,  7, -1,  7, 3,  1, -1] |
| 2 | 8 | 12 | (8, 4) | (0, 2, 4, 2) | [ 8,  0, -1, 6, 4, 1, 2, -1, 6, 2, 8, 7, -1, 7, 1, 4, -1, 1, 5, 3, 2, -1, 8, 2, 3, -1, 8, 5, 1, 7, -1, 4, 5, 8, -1, 3, 6, 7, -1, 7, 4, 3, -1, 4, 8, 3, -1, 6, 5, 4, -1, 3, 5, 6, -1] |
| 3 | 8 | 13 | (10, 3) | (1, 0, 3, 4) | [ 8,  0, -1, 8, 6, 1, 2, -1, 8, 2, 3, -1, 6, 5, 3, 1, -1, 4, 2, 1, -1, 2, 7, 3, -1, 1, 3, 7, -1, 5, 7, 2, -1, 4, 8, 3, -1, 4, 1, 7, -1, 3, 5, 4, -1, 5, 2, 4, -1, 8, 7, 5, 6, -1, 4, 7, 8, -1] |
| 4 | 9 | 11 | (4, 7) | (0, 5, 4) | [ 9,  0, -1, 8, 7, 1, 4, -1, 3, 4, 1, 2, -1, 2, 6, 7, 8, -1, 9, 1, 7, -1, 5, 2, 1, 9, -1, 2, 5, 6, -1, 6, 5, 4, 3, -1, 3, 2, 8, -1, 9, 7, 6, 3, -1, 4, 5, 9, 8, -1, 8, 9, 3, -1] |
| 5 | 9 | 11 | (4, 7) | (1, 3, 5) | [ 9,  0, -1, 9, 7, 1, 4, -1, 3, 4, 1, 2, -1, 2, 6, 7, 9, -1, 8, 1, 7, 6, -1, 5, 2, 1, 8, -1, 2, 5, 6, -1, 6, 5, 4, 3, -1, 3, 2, 9, -1, 8, 6, 3, -1, 4, 5, 8, 9,-1, 9, 8, 3, -1] |
| 6 | 9 | 11 | (1, 4, 3, 1) | (5, 5, 1) | [ 9,  0, -1, 9, 8, 1, 4, 7, -1, 4, 1, 2, 3, -1, 3, 2, 9, -1, 8, 6, 1, -1, 5, 2, 1, 6, -1, 9, 2, 5, 8, -1, 8, 5, 3, -1, 7, 6, 3, 9, -1, 6, 8, 3, -1, 4, 3, 5, -1, 7, 4, 5, 6, -1] |
| 7 | 9 | 11 | (1, 4, 3, 1) | (5, 5, 1) | [ 9,  0, -1, 4, 9, 8, 1, 6, -1, 1, 2, 3, 6, -1, 4, 3, 2, 9, -1, 5, 1, 8, -1, 7, 2, 1, 5, -1, 9, 2, 7, 8, -1, 8, 7, 3, -1, 5, 8, 3, -1, 6, 7, 5, 4, -1, 3, 4, 5, -1, 6, 3, 7, -1] |
| 8 | 9 | 11 | (6, 3, 2) | (2, 2, 4, 1) | [ 9,  0, -1, 6, 3, 4, 5, 1, -1, 6, 1, 8, 9, -1, 9, 4, 3, 2, -1, 5, 7, 2, 1, -1, 8, 1, 2, -1, 4, 9, 8, 7, 5, -1, 3, 7, 8, -1, 6, 9, 2, -1, 3, 8, 2, -1, 2, 7, 6, -1, 7, 3, 6, -1] |
| 9 | 9 | 12 | (6, 6) | (0, 3, 6) | [ 9,  0, -1, 6, 4, 1, 7, -1, 7, 9, 2, 3, -1, 4, 3, 2, 1, -1, 8, 7, 1, 5, -1, 5, 2, 9, 6, -1, 2, 5, 1, -1, 6, 9, 4, -1, 7, 3, 6, -1, 3, 8, 5, 6, -1, 3, 4, 8, -1, 4, 9, 8, -1, 9, 7, 8, -1] |
| 10 | 9 | 12 | (6, 6) | (2, 0, 6, 1) | [ 9,  0, -1, 9, 3, 1, 7, -1, 9, 7, 6, 2, -1, 2, 1, 3, 5, -1, 6, 7, 1, 4, -1, 6, 8, 2, -1, 1, 2, 8, -1, 3, 8, 6, -1, 4, 9, 2, 5, -1, 4, 1, 8, -1, 5, 3, 6, 4, -1, 9, 8, 3, -1, 4, 8, 9, -1] |
| 11 | 9 | 12 | (6, 6) | (2, 1, 4, 2) | [ 9,  0, -1, 2, 5, 7, 1, -1, 2, 1, 9, 3, -1, 8, 6, 7, 5, -1, 7, 6, 4, 1, -1, 9,  1, 4, 8, -1, 6, 8, 2, -1, 2, 8, 4, -1, 6, 2, 3, -1, 9, 8, 5, 3, -1, 3, 4, 6, -1,  5, 4, 3, -1, 2, 4, 5, -1] |
|  |  |  |  |  |  |

| 12 | 9 | 12 | (7, 4, 1) | (2, 1, 4, 2) | [ 9, 0, -1, 8, 5, 7, 1, 9, -1, 9, 1, 2, 3, -1, 5, 2, 6, 7, -1, 7, 6, 4, 1, -1, 1, 4, 2, -1, 6, 2, 8, -1, 8, 2, 4, -1, 9, 3, 6, 8, -1, 3, 2, 5, -1, 3, 4, 6, -1, 5, 4, 3, -1, 8, 4, 5, -1] |
|---|---|---|---|---|---|
|  |  |  |  |  |  |
| 13 | 10 | 10 | (0, 10) | (0, 10) | [ 10, 0, -1, 1, 4, 5, 2, -1, 2, 5, 6, 3, -1, 6, 8, 4, 3, -1, 4, 1, 7, 3, -1, 2, 9, 10, 1, -1, 10, 9, 6, 5, -1, 4, 8, 10, 5, -1, 8, 7, 1, 10, -1, 3, 7, 9, 2, -1, 7, 8, 6, 9, -1] |
|  |  |  |  |  |  |
| 14 | 10 | 10 | (1, 8, 1) | (2, 6, 2) | [ 10, 0, -1, 1, 4, 5, 2, -1, 2, 8, 9, 7, -1, 2, 7, 10, 1, -1, 10, 7, 6, 5, -1, 7, 9, 3, 6, -1, 9, 4, 1, 3, -1, 3, 1, 10, 8, -1, 10, 5, 4, 9, 8, -1, 8, 2, 3, -1, 2, 5, 6, 3, -1] |
|  |  |  |  |  |  |
| 15 | 10 | 10 | (2, 6, 2) | (1, 8, 1) | [ 10, 0, -1, 1, 4, 5, 2, -1, 2, 8, 9, 7, -1, 2, 7, 10, 1, -1, 10, 7, 6, 5, -1, 7, 9, 4, 3, 6, -1, 4, 1, 3, -1, 4, 9, 8, 10, 5, -1, 8, 3, 1, 10, -1, 8, 2, 6, 3, -1, 2, 5, 6, -1] |
|  |  |  |  |  |  |
| 16 | 10 | 10 | (2, 6, 2) | (2, 6, 2) | [ 10, 0, -1, 1, 4, 5, 2, -1, 2, 5, 6, 3, -1, 1, 2, 8, 9, -1, 9, 8, 6, 5, -1, 1, 9, 10, 3, -1, 10, 9, 5, 4, 7, -1, 3, 10, 7, 8, 2, -1, 7, 6, 8, -1, 3, 6, 1, -1, 6, 7, 4, 1, -1] |
|  |  |  |  |  |  |
| 17 | 10 | 10 | (2, 6, 2) | (2, 6, 2) | [ 10, 0, -1, 1, 4, 5, 2, -1, 2, 5, 6, 3, -1, 1, 2, 9, 10, -1, 10, 9, 6, 5, -1, 4, 8, 7, 10, 5, -1, 7, 1, 10, -1, 3, 7, 8, 9, 2, -1, 8, 6, 9, -1, 4, 1, 6, 8, -1, 1, 7, 3, 6, -1] |
|  |  |  |  |  |  |
| 18 | 10 | 10 | (2, 6, 2) | (2, 6, 2) | [ 10, 0, -1, 1, 4, 5, 2, -1, 8, 2, 5, 6, 3, -1, 2, 8, 9, 10, -1, 1, 7, 8, 3, -1, 7, 4, 9, 8, -1, 1, 2, 10, 7, -1, 10, 9, 4, 3, 6, -1, 4, 1, 3, -1, 6, 5, 7, 10, -1, 5, 4, 7, -1] |
|  |  |  |  |  |  |
| 19 | 10 | 10 | (3, 4, 3) | (3, 4, 3) | [ 10, 0, -1, 1, 4, 5, 2, -1, 2, 5, 6, 3, -1, 10, 2, 3, 8, 7, -1, 9, 5, 4, 8, -1, 2, 10, 9, 1, -1, 10, 7, 6, 5, 9, -1, 9, 8, 3, -1, 3, 1, 9, -1, 3, 6, 1, -1, 6, 7, 8, 4, 1, -1] |
|  |  |  |  |  |  |
| 20 | 10 | 10 | (3, 4, 3) | (3, 4, 3) | [ 10, 0, -1, 1, 4, 5, 2, -1, 2, 7, 10, 8, 6, -1, 6, 8, 4, 3, -1, 4, 1, 3, -1, 2, 6, 9, 1, -1, 9, 6, 5, -1, 2, 5, 3, 7, -1, 5, 6, 3, -1, 1, 9, 10, 7, 3, -1, 10, 9, 5, 4, 8, -1] |
|  |  |  |  |  |  |
| 21 | 10 | 10 | (3, 5, 1, 1) | (3, 5, 1, 1) | [ 10, 0, -1, 1, 4, 5, 2, -1, 10, 8, 3, 9, 5, 4, -1, 3, 1, 7, 9, -1, 3, 6, 1, -1, 6, 4, 1, -1, 1, 2, 8, 10, 7, -1, 10, 4, 6, 7, -1, 6, 3, 8, 2, -1, 7, 6, 9, -1, 6, 2, 5, 9, -1] |
|  |  |  |  |  |  |
| 22 | 10 | 10 | (4, 2, 4) | (4, 2, 4) | [ 10, 0, -1, 1, 4, 5, 2, -1, 9, 2, 3, 8, 7, -1, 5, 6, 2, -1, 6, 3, 2, -1, 2, 9, 10, 1, -1, 10, 9, 7, 6, 5, -1, 4, 8, 3, 10, 5, -1, 3, 1, 10, -1, 8, 4, 1, 6, 7, -1, 1, 1, 3, 6, -1] |
|  |  |  |  |  |  |
| 23 | 10 | 11 | (2, 9) | (1, 6, 3) | [ 10, 0, -1, 2, 5, 6, 3, -1, 7, 9, 1, -1, 9, 7, 6, 5, -1, 1, 4, 10, 7, -1, 10, 4, 5, 2, -1, |
| | | | | | 1, 9, 8, 3, -1, 9, 5, 4, 8, -1, 10, 2, 3, 8, -1, 8, 6, 7, 10, -1, 8, 4, 6, -1, 4, 1, 3, 6, -1] |
| 24 | 10 | 11 | (3, 7, 1) | (3, 2, 5) | [ 10, 0, -1, 9, 6, 1, 2, 10, -1, 2, 1, 8, 3, -1, 9, 7, 5, 6, -1, 6, 5, 4, 1, -1, 8, 1, 4, 7, -1, 5, 7, 2, -1, 2, 7, 4, 10, -1, 5, 2, 3, -1, 8, 7, 9, 3, -1, 3, 4, 5, -1, 10, 4, 3, 9, -1] |
| 25 | 10 | 11 | (3, 7, 1) | (3, 3, 3, 1) | [ 10, 0, -1, 9, 5, 8, 1, 10, -1, 10, 1, 2, 3, -1, 5, 2, 7, 8, -1, 8, 7, 4, 1, -1, 1, 4, 2, -1, 7, 2, 9, -1, 9, 2, 4, -1, 10, 3, 7, 9, -1, 3, 2, 5, 6, -1, 3, 6, 4, 7, -1, 9, 4, 6, 5, -1] |
| 26 | 10 | 11 | (4, 5, 2) | (0, 8, 2) | [ 10, 0, -1, 10, 1, 7, 5, -1, 1, 4, 2, 8, -1, 4, 9, 5, 2, -1, 6, 7, 1, 8, 9, -1, 9, 4, 6, -1, 10, 3, 6, 4, -1, 4, 1, 10, -1, 3, 10, 5, 9, 8, -1, 8, 2, 3, -1, 7, 6, 3, 2, -1, 2, 5, 7, -1] |
| 27 | 10 | 11 | (4, 5, 2) | (1, 6, 3) | [ 10, 0, -1, 10, 3, 6, 4, 1, -1, 10, 1, 7, 5, -1, 1, 4, 2, 8, -1, 4, 9, 5, 2, -1, 7, 1, 8, -1, 8, 9, 6, 7, -1, 9, 4, 6, -1, 3, 10, 5, 9, 8, -1, 8, 2, 3, -1, 7, 6, 3, 2, -1, 2, 5, 7, -1] |
| 28 | 10 | 11 | (4, 6, 0, 1) | (3, 3, 3, 1) | [ 10, 0, -1, 7, 10, 9, 6, 1, 8, -1, 8, 1, 2, 3, -1, 9, 2, 5, 6, -1, 6, 5, 4, 1, -1, 1, 4, 2, -1, 5, 2, 7, -1, 7, 2, 4, 10, -1, 8, 3, 5, 7, -1, 3, 2, 9, -1, 3, 4, 5, -1, 10, 4, 3, 9, -1] |
| 29 | 10 | 11 | (5, 3, 3) | (1, 6, 3) | [ 10, 0, -1, 3, 6, 4, 1, -1, 8, 2, 5, 6, 3, -1, 1, 2, 8, 10, 7, -1, 10, 6, 5, 7, -1, 8, 3, 9, -1, 6, 10, 4, -1, 10, 8, 9, 4, -1, 3, 1, 7, -1, 7, 9, 3, -1, 7, 5, 9, -1, 5, 2, 1, 4, 9, -1] |
| 30 | 10 | 11 | (5, 3, 3) | (3, 3, 3, 1) | [ 10, 0, -1, 2, 5, 6, 3, -1, 3, 1, 7, 10, -1, 6, 5, 7, 8, 9, -1, 7, 1, 8, -1, 4, 1, 6, 9, -1, 1, 3, 6, -1, 7, 5, 2, 4, 10, -1, 1, 4, 2, -1, 2, 8, 1, -1, 2, 3, 8, -1, 3, 10, 4, 9, 8, -1] |
| 31 | 10 | 11 | (5, 3, 3) | (4, 0, 6) | [ 10, 0, -1, 9, 5, 1, 8, -1, 9, 8, 2, -1, 2, 1, 5, 6, 4, -1, 8, 1, 3, -1, 7, 10, 2, 8, -1, 1, 2, 10, -1, 9, 10, 7, 6, 5, -1, 3, 9, 2, 4, -1, 3, 1, 10, -1, 4, 6, 7, 8, 3, -1, 3, 10, 9, -1] |
| 32 | 11 | 9 | (0, 5, 4) | (4, 7) | [ 11, 0, -1, 1, 4, 5, 2, -1, 8, 2, 5, 6, 3, -1, 2, 8, 11, 9, 7, -1, 2, 7, 10, 1, -1, 10, 7, 6, 5, -1, 7, 9, 3, 6, -1, 9, 4, 1, 3, -1, 1, 10, 11, 8, 3, -1, 11, 10, 5, 4, 9, -1] |
| 33 | 11 | 9 | (1, 3, 5) | (4, 7) | [ 11, 0, -1, 1, 4, 5, 2, -1, 8, 2, 5, 6, 3, -1, 2, 8, 11, 9, 7, -1, 2, 7, 10, 1, -1, 10, 7, 6, 5, -1, 7, 9, 4, 3, 6, -1, 4, 1, 3, -1, 4, 9, 11, 10, 5, -1, 11, 8, 3, 1, 10, -1] |

| 34 | 11 | 9 | (1, 4, 3, 1) | (5, 5, 1) | [ 11, 0, -1, 1, 4, 5, 2, -1, 2, 5, 6, 3, -1, 10, 2, 3, 8, 9, 7, -1, 2, 10, 11, 1, -1, 11, 10, 7, 6, 5, -1, 4, 9, 8, 11, 5, -1, 8, 1, 11, -1, 4, 1, 6, 7, 9, -1, 1, 8, 3, 6, -1] |
|----|----|---|--------------|-----------|-----|
| 35 | 11 | 9 | (1, 4, 3, 1) | (5, 5, 1) | [ 11, 0, -1, 1, 4, 5, 2, -1, 2, 5, 6, 3, -1, 11, 2, 3, 10, 8, 7, -1, 3, 1, 9, 10, -1, 10, 9, 5, 4, 8, -1, 2, 11, 9, 1, -1, 11, 7, 6, 5, 9, -1, 8, 4, 1, 6, 7, -1, 1, 3, 6, -1] |
| 36 | 11 | 9 | (2, 2, 4, 1) | (6, 3, 2) | [ 11, 0, -1, 1, 4, 5, 2, -1, 10, 8, 3, 9, 5, 4, -1, 3, 1, 11, 7, 9, -1, 3, 6, 1, -1, 6, 4, 1, -1, 1, 2, 8, 10, 11, -1, 11, 10, 4, 6, 7, -1, 7, 6, 2, 5, 9, -1, 6, 3, 8, 2, -1] |
| 37 | 11 | 10 | (0, 8, 2) | (4, 5, 2) | [ 11, 0, -1, 1, 4, 5, 2, -1, 2, 5, 6, 3, -1, 3, 1, 9, 10, -1, 10, 9, 5, 4, 8, -1, 2, 11, 9, 1, -1, 11, 7, 6, 5, 9, -1, 2, 3, 10, 8, -1, 8, 7, 11, 2, -1, 8, 4, 1, 7, -1, 1, 3, 6, 7, -1] |
| 38 | 11 | 10 | (1, 6, 3) | (2, 9) | [ 11, 0, -1, 2, 5, 6, 3, -1, 1, 7, 11, 6, 5, -1, 8, 1, 5, 9, -1, 6, 4, 3, -1, 4, 1, 8, 3, -1, 1, 4, 10, 2, 7, -1, 10, 9, 5, 2, -1, 3, 8, 11, 7, 2, -1, 8, 9, 10, 11, -1, 10, 4, 6, 11, -1] |
| 39 | 11 | 10 | (1, 6, 3) | (4, 5, 2) | [ 11, 0, -1, 1, 4, 5, 2, -1, 8, 2, 5, 6, 3, -1, 10, 2, 8, 9, 7, -1, 2, 10, 11, 1, -1, 11, 10, 7, 6, 5, -1, 4, 9, 11, 5, -1, 11, 9, 8, 3, -1, 3, 1, 11, -1, 3, 6, 7, 1, -1, 7, 9, 4, 1, -1] |
| 40 | 11 | 10 | (1, 6, 3) | (5, 3, 3) | [ 11, 0, -1, 2, 5, 6, 3, -1, 1, 7, 8, 6, 5, -1, 11, 9, 1, 5, -1, 11, 5, 4, 10, -1, 5, 2, 4, -1, 2, 7, 1, 4, -1, 8, 10, 4, 3, 6, -1, 4, 1, 9, 3, -1, 8, 7, 2, 11, 10, -1, 2, 3, 9, 11, -1] |
| 41 | 11 | 10 | (3, 2, 5) | (3, 7, 1) | [ 11, 0, -1, 1, 4, 5, 2, -1, 1, 2, 8, 6, 7, -1, 2, 9, 11, 10, 8, -1, 4, 10, 11, 7, 5, -1, 8, 10, 4, 3, 6, -1, 4, 1, 3, -1, 6, 3, 9, 5, 7, -1, 9, 2, 5, -1, 3, 1, 11, 9, -1, 1, 7, 11, -1] |
| 42 | 11 | 10 | (3, 3, 3, 1) | (3, 7, 1) | [ 11, 0, -1, 1, 4, 5, 2, -1, 9, 2, 5, 7, 6, 3, -1, 2, 9, 10, 8, -1, 4, 10, 11, 7, 5, -1, 1, 11, 9, 3, -1, 11, 10, 9, -1, 8, 10, 4, 3, 6, -1, 4, 1, 3, -1, 7, 11, 1, 8, 6, -1, 1, 2, 8, -1] |
| 43 | 11 | 10 | (3, 3, 3, 1) | (4, 6, 0, 1) | [ 11, 0, -1, 1, 4, 5, 2, -1, 9, 2, 5, 7, 6, 3, -1, 2, 9, 11, 10, 8, -1, 4, 10, 11, 7, 5, -1, 8, 10, 4, 3, 6, -1, 4, 1, 3, -1, 7, 1, 8, 6, -1, 1, 2, 8, -1, 3, 1, 11, 9, -1, 1, 7, 11, -1] |
| 44 | 11 | 10 | (3, 3, 3, 1) | (5, 3, 3) | [ 11, 0, -1, 1, 4, 5, 2, -1, 9, 2, 5, 7, 6, 3, -1, 2, 9, 10, 8, -1, 10, 4, 8, -1, 8, 4, 3, |
| | | | | | 6, -1, 4, 1, 3, -1, 3, 1, 11, 10, 9, -1, 11, 7, 5, 4, 10, -1, 7, 11, 1, 8, 6, -1, 1, 2, 8, -1] |
| 45 | 11 | 10 | (4, 0, 6) | (5, 3, 3) | [ 11, 0, -1, 3, 6, 9, 4, 1, -1, 8, 2, 5, 6, 3, -1, 9, 11, 8, 10, 4, -1, 8, 3, 10, -1, 3, 1, 7, -1, 7, 10, 3, -1, 7, 5, 10, -1, 5, 2, 1, 4, 10, -1, 2, 8, 11, 7, 1, -1, 11, 9, 6, 5, 7, -1] |
| 46 | 11 | 11 | (0, 11) | (0, 11) | [ 11, 0, -1, 3, 6, 4, 1, -1, 1, 4, 5, 2, -1, 3, 1, 7, 9, -1, 9, 10, 8, 3, -1, 10, 9, 5, 4, -1, 7, 1, 2, 8, -1, 6, 3, 8, 2, -1, 6, 11, 10, 4, -1, 11, 7, 8, 10, -1, 7, 11, 5, 9, -1, 11, 6, 2, 5, -1] |
| 47 | 12 | 8 | (0, 2, 4, 2) | (8, 4) | [ 12, 0, -1, 3, 6, 4, 1, -1, 1, 4, 5, 2, -1, 12, 8, 9, 10, 5, 4, -1, 9, 3, 1, 11, 7, 10, -1, 6, 7, 11, 12, 4, -1, 12, 11, 1, 2, 8, -1, 8, 2, 6, 3, 9, -1, 2, 5, 10, 7, 6, -1] |
| 48 | 12 | 9 | (0, 3, 6) | (6, 6) | [ 12, 0, -1, 3, 6, 4, 1, -1, 2, 5, 6, 3, -1, 8, 2, 3, 10, 9, -1, 3, 1, 7, 11, 10, -1, 6, 5, 9, 4, -1, 5, 7, 1, 8, 9, -1, 12, 2, 8, 1, 4, -1, 9, 10, 11, 12, 4, -1, 12, 11, 7, 5, 2, -1] |
| 49 | 12 | 9 | (2, 0, 6, 1) | (6, 6) | [ 12, 0, -1, 3, 6, 9, 4, 1, -1, 8, 2, 5, 6, 3, -1, 3, 1, 11, 7, 10, -1, 9, 12, 8, 10, 4, -1, 8, 3, 10, -1, 7, 5, 2, 4, 10, -1, 2, 1, 4, -1, 7, 11, 12, 9, 6, 5, -1, 12, 11, 1, 2, 8, -1] |
| 50 | 12 | 9 | (2, 1, 4, 2) | (6, 6) | [ 12, 0, -1, 1, 4, 5, 2, -1, 12, 1, 2, 8, 6, 7, -1, 2, 9, 11, 10, 8, -1, 4, 10, 11, 12, 7, 5, -1, 12, 11, 9, 3, 1, -1, 8, 10, 4, 3, 6, -1, 4, 1, 3, -1, 6, 3, 9, 5, 7, -1, 9, 2, 5, -1] |
| 51 | 12 | 9 | (2, 1, 4, 2) | (7, 4, 1) | [ 12, 0, -1, 1, 4, 5, 2, -1, 9, 2, 5, 7, 6, 3, -1, 2, 9, 11, 10, 8, -1, 4, 10, 11, 12, 7, 5, -1, 12, 11, 9, 3, 1, -1, 8, 10, 4, 3, 6, -1, 4, 1, 3, -1, 7, 12, 1, 8, 6, -1, 1, 2, 8, -1] |
| 52 | 13 | 8 | (1, 0, 3, 4) | (10, 3) | [ 13, 0, -1, 3, 6, 9, 4, 1, -1, 1, 4, 11, 7, 5, 2, -1, 10, 8, 2, 5, 6, 3, -1, 13, 8, 10, 11, 4, 9, -1, 11, 10, 3, 12, 7, -1, 12, 3, 1, -1, 6, 5, 7, 12, 13, 9, -1, 13, 12, 1, 2, 8, -1] |
| 53 | 14 | 7 | (0, 0, 0, 7) | (14, 0) | [ 14, 0, -1, 3, 6, 10, 12, 4, 1, -1, 9, 7, 1, 4, 5, 2, -1, 11, 2, 5, 8, 6, 3, -1, 13, 7, 9, 10, 6, 8, -1, 9, 2, 11, 14, 12, 10, -1, 1, 7, 13, 14, 11, 3, -1, 14, 13, 8, 5, 4, 12, -1] |

Table 17 (Pages 30-34) Partial set of Diminimal Toroidal Polyhedral Maps

| Size of Set | Diminimal TPMs Produced | New Diminimal TPMs produced |
|---|---|---|
| 5 | 1 | 0 |
| 5 | 0 | 0 |
| 5 | 2 | 0 |
| 20 | 0 | 0 |
| 20 | 0 | 0 |
| 20 | 3 | 0 |
| 40 | 0 | 0 |
| 40 | 2 | 1 |
| 40 | 1 | 0 |

Table 18 Description of subsets processed to completion.

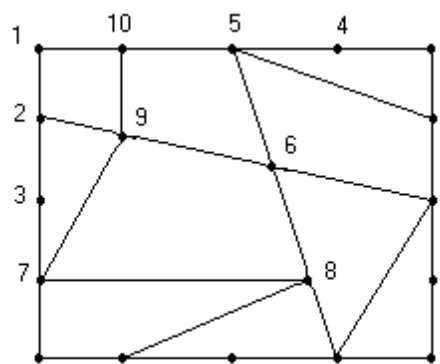

TPM 13

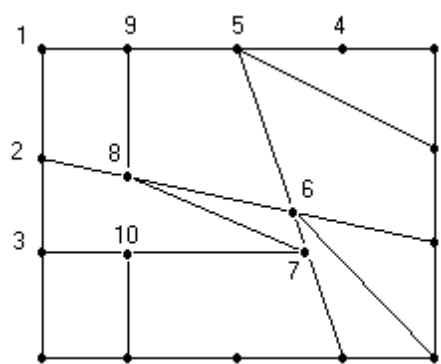

TPM 16

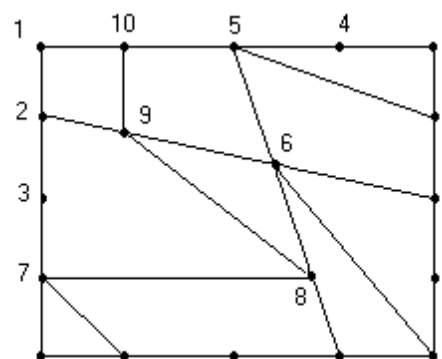

TPM 17

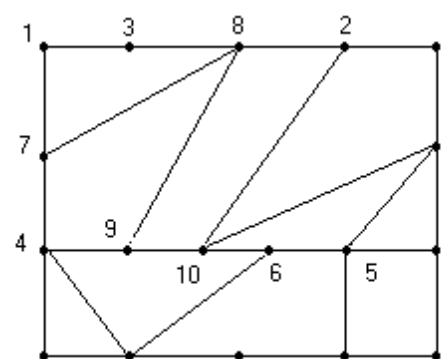

TPM 18

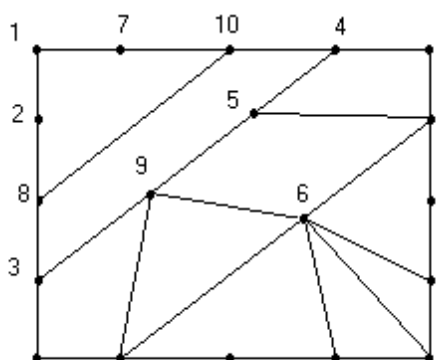

TPM 21

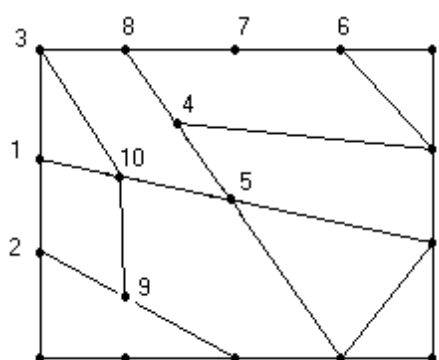

TPM 22

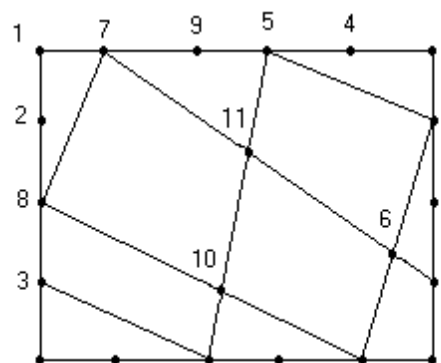

TPM 46

Figure  19  Each of the pictured TPMs is self dual.
The number under each map is its number in Table  17.

| Map # | Dual to Map # | Map # | Dual to Map # |
|-------|---------------|-------|---------------|
| 1 | 53 | 17 | self dual |
| 2 | 47 | 18 | self dual |
| 3 | 52 | 19 | 20 |
| 4 | 32 | 21 | self dual |
| 5 | 33 | 22 | self dual |
| 6 | 34 | 23 | 38 |
| 7 | 35 | 24 | 41 |
| 8 | 36 | 25 | 42 |
| 9 | 48 | 26 | 37 |
| 10 | 49 | 27 | 39 |
| 11 | 50 | 28 | 43 |
| 12 | 51 | 29 | 40 |
| 13 | self dual | 30 | 44 |
| 14 | 15 | 31 | 45 |
| 16 | self dual | 46 | self dual |

Table  20  Dual Pairs and Self Dual Maps (Numbers refer to those numbers in Table  17).

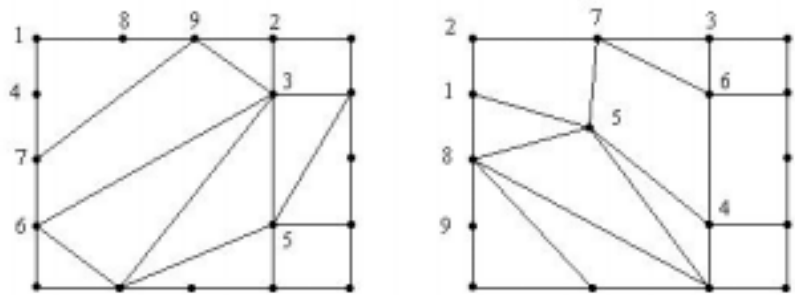

Figure  21  The two pictured TPMs have isomorphic graphs.
The isomorphism carries vertices  1  through  9  in the
left hand map to vertices  1,  6,  3,  2,  7,  8,  9,  5,  and  4
respectively in the right hand map.

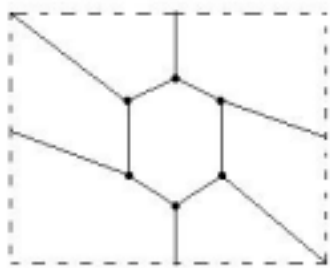

Figure  22  An embedding of $K_{(3,3)}$  in the torus.

| Map # | Dual to Map # | Map # | Dual to Map # |
|-------|---------------|-------|---------------|
| 1 | 53 | 17 | self dual |
| 2 | 47 | 18 | self dual |
| 3 | 52 | 19 | 20 |
| 4 | 32 | 21 | self dual |
| 5 | 33 | 22 | self dual |
| 6 | 34 | 23 | 38 |
| 7 | 35 | 24 | 41 |
| 8 | 36 | 25 | 42 |
| 9 | 48 | 26 | 37 |
| 10 | 49 | 27 | 39 |
| 11 | 50 | 28 | 43 |
| 12 | 51 | 29 | 40 |
| 13 | self dual | 30 | 44 |
| 14 | 15 | 31 | 45 |
| 16 | self dual | 46 | self dual |

Table  20  Dual Pairs and Self Dual Maps (Numbers refer to those numbers in Table  17).